\documentclass[11pt,reqno,tbtags]{amsart}
\usepackage[numbers]{natbib}
\usepackage{tabls} 
\usepackage{amssymb}
\numberwithin{equation}{section}

\newtheorem*{property*}{Property \csname @currentlabel\endcsname}

\newtheorem{theorem}{Theorem}[section]
\newtheorem{lemma}[theorem]{Lemma}

\newtheorem{corollary}[theorem]{Corollary}

\theoremstyle{definition}
\newtheorem{example}[theorem]{Example}
\newtheorem*{definition}{Definition}

\newtheorem{remark}{Remark}[section]

\newtheorem*{acks}{Acknowledgements}

\theoremstyle{remark}

\newenvironment{romenumerate}{\begin{enumerate}
 }{\end{enumerate}}

\newcounter{oldenumi}
{\setcounter{oldenumi}{\value{enumi}}
\begin{romenumerate} \setcounter{enumi}{\value{oldenumi}}}
{\end{romenumerate}}

\newcounter{thmenumerate}
\newenvironment{thmenumerate}
{\setcounter{thmenumerate}{0}%
 \def\item{\par
 \refstepcounter{thmenumerate}\textup{(\roman{thmenumerate})\enspace}}
}
{}

\newcounter{xenumerate}   

\newcommand{\refT}[1]{Theorem~\ref{#1}}
\newcommand{\refC}[1]{Corollary~\ref{#1}}
\newcommand{\refL}[1]{Lemma~\ref{#1}}
\newcommand{\refR}[1]{Remark~\ref{#1}}
\newcommand{\refS}[1]{Section~\ref{#1}}
\newcommand{\refand}[2]{\ref{#1} and~\ref{#2}}

\newcommand\marginal[1]{\marginpar{\raggedright\parindent=0pt\tiny #1}}

\begingroup
  \divide\count255 by 60
  \multiply\count255 by -60
  \advance\count255 by \time
  \ifnum \count255 < 10 \xdef\klockan{\the\count1.0\the\count255}
  \else\xdef\klockan{\the\count1.\the\count255}\fi
\endgroup

\newcommand{\sumi}{\sum_{i=0}^\infty}

\newcommand\set[1]{\ensuremath{\{#1\}}}

\newcommand\xpar[1]{(#1)}
\newcommand\bigpar[1]{\bigl(#1\bigr)}
\newcommand\Bigpar[1]{\Bigl(#1\Bigr)}

\newcommand\lrpar[1]{\left(#1\right)}

\newcommand\bigabs[1]{\bigl|#1\bigr|}

\def\rompar(#1){\textup(#1\textup)}    

\def\xexp(#1){e^{#1}}

\newcommand\ntoo{\ensuremath{{n\to\infty}}}

\newcommand\iid{i.i.d.\spacefactor=1000}    
\newcommand\ie{i.e.\spacefactor=1000}

\newcommand\cf{cf.\spacefactor=1000}
\newcommand{\as}{a.s.\spacefactor=1000}
\newcommand{\aex}{a.e.\spacefactor=1000}

\newcommand{\tend}{\longrightarrow}
\newcommand\dto{\overset{\mathrm{d}}{\tend}}
\newcommand\pto{\overset{\mathrm{p}}{\tend}}
\newcommand\asto{\overset{\mathrm{a.s.}}{\tend}}
\newcommand\eqd{\overset{\mathrm{d}}{=}}

\newcommand\bbN{\mathbb N}

\newcounter{CC} 
\newcounter{cc}

\newcommand\E{\operatorname{\mathbb E{}}}
\renewcommand\P{\operatorname{\mathbb P{}}}

\newcommand\Be{\operatorname{Be}}

\newcommand\ga{\alpha}
\newcommand\gb{\beta}
\newcommand\gd{\delta}

\newcommand\gf{\varphi}

\newcommand\gG{\Gamma}

\newcommand\gs{\sigma}

\newcommand\cB{\mathcal B}
\newcommand\cD{\mathcal D}

\newcommand\cF{\mathcal F}
\newcommand\cL{{\mathcal L}}
\newcommand\cP{\mathcal P}
\newcommand\cS{{\mathcal S}}
\newcommand\cU{{\mathcal U}}
\newcommand\cW{\mathcal W}

\newcommand\ett[1]{\boldsymbol1[#1]} 
\newcommand\bigett[1]{\boldsymbol1\bigl[#1\bigr]} 

\def\[#1]{[\![#1]\!]}

\newcommand\qw{^{-1}}

\renewcommand{\=}{:=}

\newcommand\oi{[0,1]}

\newcommand\dd{\,\textup{d}}

\newcommand\lhs{left hand side}
\newcommand\rhs{right hand side}

\newcommand{\cupn}{\bigcup_{n=1}^\infty}
\newcommand{\xk}{[k]}
\newcommand{\vvk}{v_1,\dots,v_k}
\newcommand{\vvki}{v'_1,\dots,v'_k}
\newcommand{\Lovasz}{Lov\'asz}
\newcommand{\Lovaszetal}{Lov\'asz and Szegedy}
\newcommand{\Borgsetal}{Borgs, Chayes, Lov\'asz, S\'os and Vesztergombi}
\newcommand{\tinj}{t_{\mathrm{inj}}}
\newcommand{\tind}{t_{\mathrm{ind}}}
\newcommand{\cux}{\cU^*}
\newcommand{\cuxq}{\overline{\cux}}
\newcommand{\cuq}{\overline{\cU}}
\newcommand{\cuu}{\cU^+}
\newcommand{\cuoo}{\cU_\infty}
\newcommand{\cloo}{\cL_\infty}
\newcommand{\oiu}{\oi^{\cU}}
\newcommand{\oiuu}{\oi^{\cuu}}
\newcommand{\Koo}{K_\infty}
\newcommand{\tauu}{\tau^+}
\newcommand{\tauinj}{\tau_{\mathrm{inj}}}
\newcommand{\tauind}{\tau_{\mathrm{ind}}}
\newcommand{\dcut}{\gd_\square}
\newcommand\hG{\widehat G}
\newcommand\hgn{\widehat {G_n}}
\newcommand{\rest}[1]{|_{[#1]}}

\newcommand{\restx}[1]{|_{#1}}
\newcommand{\restkk}{|_{[k_1]\times[k_2]}}
\newcommand{\restnn}{|_{[n_1]\times[n_2]}}
\newcommand{\exch}{exchangeable}
\newcommand{\permgs}{\circ\gs}
\newcommand{\Hgs}{H\permgs}
\newcommand{\BU}{\cB}
\newcommand{\BUxx}[2]{\cB_{#1#2}}
\newcommand{\BUnn}{\BUxx{n_1}{n_2}}
\newcommand{\BL}{\cB^L}
\newcommand{\BLxx}[2]{\cB^L_{#1#2}}
\newcommand{\BLnn}{\BLxx{n_1}{n_2}}
\newcommand{\BLkk}{\BLxx{k_1}{k_2}}
\newcommand{\cbx}{\cB^*}
\newcommand{\cbxq}{\overline{\cbx}}
\newcommand{\cbq}{\overline{\cB}}
\newcommand{\cbb}{\cB^+}
\newcommand{\cboo}{\BUxx{\infty}{\infty}}
\newcommand{\cbloo}{\BLxx{\infty}{\infty}}
\newcommand{\oib}{\oi^{\cB}}
\newcommand{\oibb}{\oi^{\cbb}}
\newcommand{\DU}{\cD}
\newcommand{\DUx}[1]{\cD_{#1}}
\newcommand{\DUn}{\DUx{n}}
\newcommand{\DL}{\cD^L}
\newcommand{\DLx}[1]{\cD^L_{#1}}
\newcommand{\DLn}{\DLx{n}}

\newcommand{\cdx}{\cD^*}
\newcommand{\cdxq}{\overline{\cdx}}
\newcommand{\cdq}{\overline{\cD}}
\newcommand{\cdd}{\cD^+}
\newcommand{\cdoo}{\DUx{\infty}}

\newcommand{\oid}{\oi^{\cD}}
\newcommand{\oidd}{\oi^{\cdd}}
\newcommand{\cw}{\cW}
\newcommand{\cws}{\cW_{\mathsf s}}
\newcommand{\cwsq}{\widehat\cW_{\mathsf s}}
\newcommand{\gwx}[1]{G(#1,W)}
\newcommand{\gwoo}{\gwx\infty}
\newcommand{\gwn}{\gwx{n}}
\newcommand{\gwxx}[2]{G(#1,#2,W)}
\newcommand{\gwoooo}{\gwxx{\infty}{\infty}}
\newcommand{\gwnn}{\gwxx{n_1}{n_2}}
\newcommand{\gbwx}[1]{G(#1,\bW)}
\newcommand{\gbwoo}{\gbwx\infty}
\newcommand{\gbwn}{\gbwx{n}}
\newcommand{\gbwpx}[1]{G(#1,\bW,p)}
\newcommand{\gbwpoo}{\gbwpx\infty}
\newcommand{\gbwpn}{\gbwpx{n}}
\newcommand{\uw}{\Gamma_W}
\newcommand{\uww}{\Gamma''_W}
\newcommand{\ubw}{\Gamma_{\bW}}
\newcommand{\ubwp}{\Gamma_{\bW,p}}
\newcommand{\sn}{\mathfrak S_n}
\newcommand{\W}[1]{W\circ#1}
\newcommand{\xio}{\xi_\emptyset}
\newcommand{\xii}{\xi_i}
\newcommand{\xij}{\xi_j}
\newcommand{\xiij}{\xi_{ij}}
\newcommand{\uoi}{U(0,1)}
\newcommand{\gko}{g_{k,0}}
\newcommand{\gki}{g_{k,1}}
\newcommand{\gkii}{g_{k,2}}
\newcommand{\gii}{g_{1,1}}
\newcommand{\gji}{g_{2,1}}
\newcommand{\wab}{W_{\ga\gb}}
\newcommand{\wo}{w}
\newcommand{\bW}{\mathbf{W}}
\newcommand{\setoi}{\set{0,1}}
\newcommand{\xx}{Y}
\newcommand{\csp}{\cS}
\newcommand{\WW}{\cW_5}
\newcommand{\WWW}{\cW_4}

\newcommand\REM[1]{{\raggedright\texttt{[#1]}\par\marginal{XXX}}}

\newcommand\urladdrx[1]{{\urladdr{\def~{{\tiny$\sim$}}#1}}}

\begin{document}
\title
{Graph limits and exchangeable random graphs}

\date{December 10, 2007} 

\author{Persi Diaconis}
\address{Department of Mathematics,
Stanford University, Stanford California 94305, USA
and  
D\'epartement de Math\'ematiques,
Universit\'e de Nice - Sophia Antipolis,
Parc Valrose,
06108 Nice Cedex 02,
France}

\author{Svante Janson}
\address{Department of Mathematics, Uppsala University, PO Box 480,
SE-751~06 Uppsala, Sweden}
\email{svante.janson@math.uu.se}
\urladdrx{http://www.math.uu.se/~svante/}

\subjclass[2000]{} 

\begin{abstract} 
We develop a clear connection between deFinetti's theorem for
exchangeable arrays (work of Aldous--Hoover--Kallenberg) and the
emerging area of graph limits (work of Lov\'asz and many
coauthors). Along the way, we translate the graph theory into more
classical probability.
\end{abstract}

\maketitle

\section{Introduction}\label{S:intro}

DeFinetti's profound contributions are now woven into many parts of
probability, statistics and philosophy. Here we show how developments
from deFinetti's work on partial exchangeability have a direct link to
the recent development of a limiting theory for large graphs. This
introduction first recalls the theory of exchangeable arrays
(\refS{secA}). Then, the subject of graph limits is outlined
(\refS{secB}). Finally, the link between these ideas, which forms the
bulk of this paper, is outlined (\refS{secC}).

\subsection{Exchangeability, partial exchangeability and exchangeable
  arrays}\label{secA}

Let $\{ X_i \}$, $1 \leq i < \infty$, be a sequence of binary random
variables. They are {\em exchangeable} if
\begin{equation*}
\P(X_1=e_1, \dots, X_n=e_n) = \P(X_1=e_{\sigma(1)}, \dots,
X_n=e_{\sigma(n)})
\end{equation*}
for all $n$, permutations $\sigma\in\sn$ 
and all $e_i \in \{ 0,1 \}$. The
celebrated representation theorem says
\begin{theorem}[deFinetti] \label{deF}
If\/ $\{ X_i \}$,  $1 \leq i<\infty$, is a binary exchangeable sequence, then:
\begin{thmenumerate}
\item With probability $1$, $X_\infty = \lim \frac{1}{n} (X_1 + \cdots
  + X_n)$ exists.
\item If $\mu(A)=P \{ X_\infty \in A \}$, then for all $n$ and  $e_i$,
  $1 \leq i \leq n$,
\begin{equation}
\P(X_1 = e_1, \dots, X_n=e_n) = \int_0^1 x^s (1-x)^{n-s} \mu(dx)
\label{eqP1}
\end{equation}
for $s=e_1 + \dots + e_n$.
\end{thmenumerate}
\end{theorem}

It is natural to refine and extend deFinetti's theorem to allow more
general observables ($X_i$ with values in a Polish space) and other
notions of symmetry (partial exchangeability). A definitive treatment
of these developments is given in \citet{Kallenberg:exch}. 
Of interest here is the
extension of deFinetti's theorem to two-dimensional arrays.
\begin{definition}
Let $\{ X_{ij} \}$, $1\leq i, j < \infty$, be binary random
variables. They are {\em separately exchangeable} if
\begin{equation}\label{exch}
\P(X_{ij} = e_{ij},\, 1 \leq i,j \leq n) = \P(X_{ij} =
e_{\sigma(i)\tau(j)},\, 1 \leq i,j \leq n)
\end{equation}
for all $n$, all permutations $\sigma, \tau\in\sn$ and all $e_{ij} \in \{
0,1 \}$. 
They are \emph{(jointly) exchangeable} if \eqref{exch} holds in the
special case $\tau=\gs$.

Equivalently, the array \set{X_{ij}} is jointly \exch{} if 
the array $\set{X_{\gs(i)\gs(j)}}$ has the same distribution as
$\set{X_{ij}}$ for every permutation $\gs$ of $\bbN$, and similarly
for separate exchangeability.
\end{definition}

The question of two-dimensional versions of deFinetti's theorem under
(separate) exchangeability arose from the statistical problems of two-way
analysis of variance. Early workers expected a version of (\ref{eqP1})
with perhaps a two-dimensional integral. The probabilist David Aldous
\cite{Aldous} and the logician Douglas Hoover \cite{Hoover} found that
the answer is more 
complicated.

Define a random binary array $\{ X_{ij} \}$ as follows: Let $U_i, V_j$,
$1 \leq i,j < \infty$, be independent and uniform in $[0,1]$. Let
$W(x,y)$ be a function from $[0,1]^2$ to $[0,1]$. Let $X_{ij}$ be $1$
or $0$ as a $W(U_i, V_j)$-coin comes up heads or tails. Let $P_W$ be
the probability distribution of $\{ X_{ij} \}$, $1 \leq i,j <
\infty$. The family $\{ X_{ij} \}$ is separately exchangeable because of
the symmetry of the construction. The Aldous--Hoover theorem says that
any separately exchangeable binary array is a mixture of such $P_W$:
\begin{theorem}[Aldous--Hoover]
Let $X = \{ X_{ij} \}$, $1 \leq i,j < \infty$, be a separately exchangeable
binary array. Then, there is a probability $\mu$ such that
\begin{equation*}
\P \{ X \in A \} = \int P_W (A) \mu (dW).
\end{equation*}
\end{theorem}

There is a similar result for jointly exchangeable arrays. 

The uniqueness of $\mu$ resisted understanding; if $\widehat W$ is
obtained from $W$ by a measure-preserving change of each variable,
clearly the associated process $\{ \widehat{X}_{ij} \}$ has the same joint
distribution as $\{ X_{ij} \}$. Using model theory, \citet{Hoover} was able to
show that this was the only source of non-uniqueness. A `probabilist's
proof' was finally found by Kallenberg, see 
\citep[Sect.~7.6]{Kallenberg:exch} 
for details and
references.

These results hold for higher dimensional arrays with $X_{ij}$ taking
values in a Polish space with minor change
\citep[Chap.~7]{Kallenberg:exch}. The description above has not
mentioned several elegant results of the theory. In particular,
Kallenberg's `spreadable' version of the theory replaces invariance
under a group by invariance under subsequences. A variety of tail
fields may be introduced to allow characterizing when $W$ takes
values in $\{ 0,1 \}$ \citep[Sect.~4]{DF1981}. Much more
general notions of partial exchangeability are studed in \cite{DF1984}.

\subsection{Graph limits}\label{secB}

Large graphs, both random and deterministic, abound in
applications. They arise from the internet, social networks, gene
regulation, ecology and in mathematics. It is natural to seek an
approximation theory: What does it mean for a sequence of graphs to
converge? When can a large complex graph be approximated by a small
graph?

In a sequence of papers
\cite{BCL3,BCL2,BCL1,BCL1ii,Freed,LL2006,LL2007,LSos,LSzcont,LSzszem,LSztest,LSz},
Laszlo \Lovasz{} with coauthors (listed here in order
of frequency) V.~T. S\'os, B. Szegedy, C. Borgs, J. Chayes, K. Vesztergombi,
A. Schrijver, M. Freedman
have developed a beautiful, unifying
limit theory. This sheds light on topics such as graph
homomorphisms, Szemeredi's regularity lemma, quasi-random graphs,
graph testing and extremal graph theory. Their theory has been
developed for dense graphs (number of edges comparable with the square
of number of vertices) but parallel theories for sparse graphs are
beginning to emerge \cite{BR}.

Roughly, a growing sequence of finite graphs $G_n$ converges if, for
any fixed graph $F$, the proportion of copies of $F$ in $G_n$
converges. \refS{Sdef} below has precise definitions.

\begin{example}

Define a probability distribution on graphs on $n$-vertices as
follows. Flip a $\theta$-coin for each vertex (dividing vertices into
`boys' and `girls'). Connect two boys with probability $p$. Connect
two girls with probability $p'$. Connect a boy and a girl with
probability $p''$. Thus, if $p=p'=0$, $p''=1$, we have a random bipartite
graph. If $p=p'=1$, $p''=0$, we have two disjoint complete graphs. If
$p=p'=p''$, we have the Erd\"os--Renyi model. As $n$ grows, these models
generate a sequence of random graphs which converge almost surely to a
limiting object described below.

More substantial   
examples involving random threshold graphs are in \cite{threshold}. 

\end{example}

If a sequence of graphs converges, what does it converge to? For
exchangeable random graphs (defined below), there is a limiting object
which may be thought of as a probability measure on infinite random
graphs. Suppose $W(x,y) = W(y,x)$ is a function from $[0,1]^2 \to
[0,1]$. Choose $\{U_i \}$, $1 \leq i < \infty$, independent uniformly
distributed random variables on $[0,1]$. Form an infinite random graph
by putting an edge from $i$ to $j$ with probability $W(U_i,
U_j)$. This measure on graphs (or alternatively $W$) is the limiting
object.

For the ``boys and girls'' example above, $W$ may be pictured as
\medskip
\begin{align*}
\theta&\ \begin{tabular}{|c|c|}
\hline[5pt]
 $p$ &$p''$ \\
[5pt]\hline[5pt]
 $p''$ &$p'$ \\
[5pt]\hline
\end{tabular}\\ 
0& \hskip 22pt 
{\theta}
\hskip 20pt 
1
\end{align*}

The theory developed shows that various properties of $G_n$ can be
well approximated by calculations with the limiting object. There is an
elegant characterization of these `continuous graph properties' with
applications to algorithms for graph testing (Does this graph
contain an Eulerian cycle?) or parameter estimation (What is an
approximation to the size of the maximum cut?). There is a practical
way to find useful approximations to a large graph by graphs of fixed
size \cite{BCL3}. 
This paper also contains a useful review of the current state of the
theory with 
proofs and references.

We have sketched the theory for unweighted graphs. There are
generalizations to graphs with weights on vertices and edges, to
bipartite, directed and hypergraphs. The sketch leaves out many nice
developments. For example, the useful cut metric between graphs \cite{LSz}
and connections to statistical physics \cite{BCL1ii}.

\subsection{Overview of the present paper}\label{secC}

There is an apparent similarity between the measure $P_W$ of the
Aldous--Hoover theorem and the limiting object $W$ from graph
limits. Roughly, working with symmetric $W$ gives the graph limit
theory; working with general $W$ gives directed graphs. The main
results of this paper make these connections precise.

Basic definitions are in \refS{Sdef} which introduces a probabilist's
version of graph convergence equivalent to the definition using graph
homomorphisms. \refS{Sconv} uses the well-established theory of weak
convergence of a sequence of probability measures on a metric space to
get properties of graph convergence. \refS{Sinfinite} carries things
over to infinite graphs.

The main results appear in \refS{Sexch}. This introduces exchangeable
random graphs and gives a one-to-one correspondence between infinite
exchangeable random graphs and distributions on the space of proper
graph limits (\refT{TE}), which specializes to a one-to-one 
correspondence between proper graph limits and extreme points in the
set of distributions of \exch{} random graphs (\refC{CE}).

A useful characterization of the extreme points of the
set of exchangeable random graphs is in \refT{TE2}. These results are
translated to the equivalence between proper graph limits and the
Aldous--Hoover theory in \refS{Srep}.
The non-uniqueness of the representing $W$, 
for \exch{} random graphs and for graph limits, 
is discussed in \refS{Sunique}.

The equivalence involves symmetric $W(x,y)$ and a single permutation
$\sigma$ taking $W(U_i, U_j)$ to $W(U_{\sigma(i)},
U_{\sigma(j)})$. The original Aldous--Hoover theorem, with perhaps
non-symmetric $W(x,y)$ and $W(U_i, V_j)$ to
$W(U_{\sigma(i)},V_{\tau(j)})$ 
translates to a limit theorem for bipartite
graphs. This is developed in \refS{Sbip}. 
The third case of the Aldous--Hoover theory for 
two-dimensional arrays, perhaps 
non-symmetric $W(x,y)$ and a single permutation $\gs$,
corresponds to directed graphs; this is sketched in \refS{Sdir}. 

The extensions to weighted
graphs are covered by allowing $X_{ij}$ to take general values in the
Aldous--Hoover theory. The extension to hypergraphs follows from the
Aldous--Hoover theory for higher-dimensional arrays.
(The details of these extensions are left to the reader.) 

Despite these parallels, the theories have much to contribute to each
other. The algorithmic, graph testing, Szemeredi partitioning
perspective is new to exchangeability theory. Indeed, the ``boys and
girls'' random graph was introduced to study the psychology of vision
in Diaconis--Freedman (1981). As far as we know, its graph theoretic
properties have not been studied. The various developments around
shell-fields in exchangeability, which characterize zero/one $W(x,y)$,
have yet to be translated into graph-theoretic terms.

\begin{acks}
This lecture is an extended version of a talk presented by PD at the 100th
anniversary of deFinetti's birth in Rome, 2006. We thank the
organizers. This work was partially funded by the French ANR's Chaire
d'excellence grant to PD.

SJ thanks 
Christian Borgs and Jennifer Chayes 
for inspiration from
lectures and discussions 
during the Oberwolfach meeting
`Combinatorics, Probability and Computing', held in November, 2006.
Parts of the research were completed during a visit by SJ
to the Universit\'e de Nice - Sophia Antipolis in January 2007.
\end{acks}

\section{Definitions and basic properties}\label{Sdef}

All graphs will be simple.
Infinite graphs will be important in later sections, but will always
be clearly stated to be infinite;  otherwise, graphs
will be finite.
We denote the vertex and edge sets of a graph $G$ by $V(G)$ and
$E(G)$, and the numbers of vertices and edges by $v(G)\=|V(G)|$ 
and 
$e(G)\=|E(G)|$. We consider both labelled and unlabelled graphs; the
labels will 
be the integers
$1,\dots,n$, where $n$ is the number of vertices in the graph.
A labelled graph is thus a graph with vertex set
$[n]\=\set{1,\dots,n}$ for some $n\ge1$; we let $\cL_n$ denote the set
of the $2^{\binom n2}$ labelled graphs on $[n]$ and let
$\cL\=\cupn\cL_n$.
An unlabelled graph can be regarded as a labelled graph where we
ignore the labels; formally, we define $\cU_n$, the set of unlabelled
graphs of order $n$, as the quotient set $\cL_n/\cong$ of labelled
graphs modulo isomorphisms.
We let $\cU\=\cupn\cU_n=\cL/\cong$, the set of all unlabelled graphs.

Note that we can, and often will, regard a labelled graph as an
unlabelled graph.

If $G$ is an (unlabelled) graph and $v_1,\dots,v_k$ is a sequence of
vertices in $G$, then $G(v_1,\dots,v_k)$ denotes the labelled graph
with vertex set $[k]$ where we put an edge between $i$ and $j$ if
$v_i$ and $v_j$ are adjacent in $G$. We allow the possibility that
$v_i=v_j$ for some $i$ and $j$. (In this case, there is no edge $ij$
because there are no loops in $G$.)

We let $G\xk$, for $k\ge1$, be the random graph $G(v_1,\dots,v_k)$
obtained by sampling $\vvk$ uniformly at random among the vertices of
$G$, with replacement. In other words, $\vvk$ are independent
uniformly distributed random vertices of $G$.

For $k\le v(G)$, we further let $G\xk'$ be the random graph $G(\vvki)$
where we sample $\vvki$ uniformly at random without replacement; the
sequence $\vvki$ is thus a uniformly distributed random sequence of
$k$ distinct vertices.

The graph limit theory in \cite{LSz} and subsequent papers is based on 
the study of the functional $t(F,G)$ 
which is defined for two graphs $F$ and $G$ as the proportion of all
mappings $V(F)\to V(G)$ that are
graph homomorphisms $F\to G$, \ie, map adjacent vertices to adjacent
vertices.
In probabilistic terms, $t(F,G)$ is the probability that a uniform random 
mapping $V(F)\to V(G)$ is a graph homomorphism.
Using the notation introduced above, we can, equivalently, write this
as, 
assuming that $F$ is labelled and $k=v(F)$,
\begin{equation}
  \label{t}
t(F,G)\=\P\bigpar{F\subseteq G\xk}.
\end{equation}
Note that both $F$ and $G\xk$ are graphs on $\xk$,
so the relation $F\subseteq G\xk$ is well-defined as containment of
labelled graphs on the same vertex set, \ie{} as $E(F)\subseteq E(G\xk)$.
Although the relation $F\subseteq G\xk$ may depend on the
labelling of $F$, the probability in \eqref{t} does not, by symmetry,
so $t(F,G)$ is really well defined by \eqref{t} for unlabelled $F$ and
$G$.

With $F$, $G$ and $k$ as in \eqref{t}, 
we further define, again 
following \cite{LSz} (and the notation of \cite{BCL1})
but stating the definitions in different but equivalent forms,
\begin{align}
  \tinj(F,G)&\=\P\bigpar{F\subseteq G\xk'}\label{tinj}
\intertext{and}
  \tind(F,G)&\=\P\bigpar{F= G\xk'},  \label{tind}
\end{align}
provided $F$ and $G$ are (unlabelled) graphs with $v(F)\le v(G)$.
If $v(F)>v(G)$ we set $\tinj(F,G)\=\tind(F,G)\=0$.

Since the probability that a random sample $\vvk$ of vertices in $G$
contains some repeated vertex is $\le k^2/(2v(G))$, it follows that
\cite{LSz}
\begin{equation}
  \label{a4}
\bigabs{t(F,G)-\tinj(F,G)}
\le\frac{v(F)^2}{2v(G)}.
\end{equation}
Hence, when considering asymptotics with $v(G)\to \infty$, it does
not matter whether we use $t$ or $\tinj$. Moreover, if $F\in\cL_k$,
then, as pointed out in \cite{BCL1} and \cite{LSz},
\begin{align}
  \tinj(F,G)&=\sum_{F'\in\cL_k,\;F'\supseteq F} \tind(F,G)
\label{a4a}
\intertext{and, by inclusion-exclusion,}
  \tind(F,G)&=\sum_{F'\in\cL_k,\;F'\supseteq F} (-1)^{e(F')-e(F)}\tinj(F,G).
\label{a4b}
\end{align}
Hence, the two families
$\set{\tinj(F,\cdot)}_{F\in\cU}$ and 
$\set{\tind(F,\cdot)}_{F\in\cU}$ of graph functionals
contain the same information and can replace each other.

The basic definition of \Lovasz{} and Szegedy 
\cite{LSz} and 
\Borgsetal{} \cite{BCL1} 
is that a sequence $(G_n)$ of graphs
converges if $t(F,G_n)$ converges for every graph $F$. We can express
this by considering the map $\tau:\cU\to\oiu$ defined by
\begin{equation}\label{tau}
  \tau(G)\=(t(F,G))_{F\in\cU}\in\oiu. 
\end{equation}
Then $(G_n)$ converges if and only if $\tau(G_n)$ converges in $\oiu$,
equipped with the usual product topology. Note that 
$\oiu$ is a compact metric space; as is well known, a metric can be
defined by, for example,
\begin{equation}
  \label{a5}
d\bigpar{(x_F),(y_F)} \= \sumi 2^{-i}|x_{F_i}-y_{F_i}|, 
\end{equation}
where
$F_1,F_2,\dots$ is some enumeration of all unlabelled graphs.

We define $\cux\=\tau(\cU)\subseteq\oiu$ to be the image of $\cU$
under this mapping $\tau$, and let $\cuxq$ be the closure of $\cux$ in
$\oiu$.
Thus $\cuxq$ is a compact metric space.
(For explicit descriptions of the subset $\cuxq$ of $\oiu$ as a set of
graph functionals, see
\Lovaszetal{} \cite{LSz}.)

As pointed out in \cite{LSz} and \cite{BCL1} (in equivalent terminology),
$\tau$ is not injective; for example, $\tau(K_{n,n})$ is the same for
all complete bipartite graphs $K_{n,n}$. 
Nevertheless,
as in \cite{LSz} and \cite{BCL1}, we can
consider a graph $G$ as an element of $\cux$ by identifying $G$ and
$\tau(G)$ (thus identifying graphs with the same $\tau(G)$), and then
convergence of $(G_n)$ as defined above is equivalent to convergence
in $\cuxq$. The limit is thus an element of $\cuxq$, but typically not
a graph in $\cux$.
The main result of \Lovaszetal{} \cite{LSz} is a representation of
the elements in $\cuxq$ to which we will return in \refS{Srep}.

\begin{remark}\label{Rmetric1}
  As said above, $\cuxq$ is a compact metric space, and it can be
  given several equivalent metrics.
One metric is the metric \eqref{a5} inherited from $\oiu$, which for 
graphs becomes 
$d(G,G')=\sum_i 2^{-i}|t(F_i,G)-t(F_i,G')|$.
Another metric, shown by
\Borgsetal{} \cite{BCL1} to be equivalent, 
is the cut-distance $\dcut$, see 
\cite{BCL1} for definitions.
Further characterizations of convergence of sequences of graphs in
  $\cuq$ are given in \citet{BCL1,BCL1ii}. 
\end{remark}

The identification of graphs with the same image in $\cux$
(\ie, with the same $t(F,\cdot)$ for all $F$) is sometimes elegant but
at other times inconvenient. It can be avoided if we instead let
$\cuu$ be the union of $\cU$ and some one-point set \set{*} and
consider the mapping $\tauu:\cU\to\oiuu=\oiu\times\oi$ defined by
\begin{equation}\label{tau+}
  \tauu(G)=\bigpar{\tau(G),\,v(G)\qw}.
\end{equation}
Then $\tauu$ is injective, because if $\tau(G_1)=\tau(G_2)$ for two
graphs $G_1$ and $G_2$ with the same number of vertices, then $G_1$
and $G_2$ are isomorphic and thus $G_1=G_2$ as unlabelled graphs.
(This can easily be shown directly: it follows from \eqref{t} that
$G_1[k]\eqd G_2[k]$ for every $k$, which implies $G_1[k]'\eqd G_2[k]'$
for every $k\le v(G_1)=v(G_2)$; now take $k=v(G_1)$.
It is also a consequence of
\cite[Theorem 2.7 and Theorem 2.3 or Lemma 5.1]{BCL1}.)

Consequently, we can identify $\cU$ with its image $\tauu(\cU)\subseteq\oiuu$
and define $\cuq\subseteq\oiuu$ as its closure.
It is easily seen that a sequence $(G_n)$ of graphs converges in
$\cuq$ if and only if either $v(G_n)\to\infty$ and $(G_n)$ converges
in $\cuxq$, or the sequence $(G_n)$ is constant from some $n_0$ on.
Hence, convergence in $\cuq$ is essentially the same as the
convergence considered by by \Lovaszetal{} \cite{LSz}, but without any
identification of non-isomorphic graphs of different orders. 

Alternatively, we can consider $\tauinj$ or $\tauind$ defined by
\begin{align*}
  \tauinj(G)&\=(\tinj(F,G))_{F\in\cU}\in\oiu, 
\\
  \tauind(G)&\=(\tind(F,G))_{F\in\cU}\in\oiu.
\end{align*}
It is easy to see that both $\tauinj$ and $\tauind$ are injective
mappings $\cU\to\oiu$.
(If $\tauinj(F,G_1)=\tauinj(F,G_2)$ for all $F$, we take $F=G_1$ and
$F=G_2$
and conclude $G_1=G_2$,
using our special definition above when $v(F)>v(G)$.)
Hence, we can again identify $\cU$ with its image
and consider its closure $\cuq$ in $\oiu$. Moreover, using \eqref{a4},
\eqref{a4a}, and \eqref{a4b}, it is easily shown that if $(G_n)$ is a
sequence of unlabelled graphs, then 
\begin{equation*}
  \tauu(G_n) \text{ converges}
\iff
  \tauind(G_n) \text{ converges}
\iff
  \tauinj(G_n) \text{ converges}.
\end{equation*}
Hence, the three compactifications $\overline{\tauu(\cU)}$, 
$\overline{\tauinj(\cU)}$, $\overline{\tauind(\cU)}$ are homeomorphic
and we can use any of them for $\cuq$.
We let $\cuoo\=\cuq\setminus\cU$; 
this is the set of all limit objects of sequences
$(G_n)$ in $\cU$ with $v(G_n)\to\infty$. (I.e., it is the set of all
proper graph limits.)

We will in the sequel prefer to use $\cuq$ rather than $\cuxq$, thus
not identifying some graphs of different orders, nor identifying finite graphs
with some limit objects in $\cuoo$.

For every fixed graph $F$, the functions $t(F,\cdot)$,
$\tinj(F,\cdot)$ and $\tind(F,\cdot)$ have unique continuous
extensions to $\cuq$, for which we use the same notation. 
We similarly extend $v(\cdot)\qw$ continuously
to $\cuq$ by defining $v(G)=\infty$
and thus $v(G)\qw=0$ for $G\in\cuoo\=\cuq\setminus\cU$.
Then \eqref{a4},
\eqref{a4a} and \eqref{a4b} hold for all $G\in\cuq$, where \eqref{a4}
means that 
\begin{equation}
  \label{a4x}
\tinj(F,G)=t(F,G), \qquad G\in\cuoo.
\end{equation}

Note that $\cuq$ is a compact metric space. Different, equivalent,
metrics are given by the embeddings $\tauu$, $\tauinj$, $\tauind$ into
$\oiuu$ and $\oiu$. Another equivalent metric is, by \refR{Rmetric1}
and the definition of $\tauu$, $\dcut(G_1,G_2)+|v(G_1)\qw-v(G_2)\qw|$.

We summarize the results above on convergence.

\begin{theorem}\label{T1}
  A sequence $(G_n)$ of graphs converges in the sense of \Lovaszetal{}
  \cite{LSz} if and only if it converges in the compact metric space
  $\cuxq$.
Moreover, if $v(G_n)\to\infty$, the sequence $(G_n)$ converges in this
  sense if and only if it converges in $\cuq$.
\end{theorem}

The projection 
$\pi:\oiuu=\oiu\times\oi\to\oiu$ maps $\tauu(G)$ to $\tau(G)$ for every
graph $G$, so by continuity it maps $\cuq$ into $\cuxq$.
For graph $G\in\cU$, $\pi(G)=\tau(G)$ is the object in $\cuxq$
corresponding to $G$ considered above, and we will in the sequel
denote this object by $\pi(G)$; recall that this projection
$\cU\to\cuxq$ is not injective. (We thus distinguish between a graph
$G$ and its ``ghost'' 
$\pi(G)$ in $\cuxq$. Recall that when graphs are
considered as elements of $\cuxq$ as in \cite{LSz} and \cite{BCL1},
certain graphs are identified with each other; we avoid this.)
On the other hand, an element $G$ of $\cuq$ is by definition
determined by $\tau(G)$ and $v(G)\qw$, 
\cf{} \eqref{tau+},
so the restriction $\pi:\cU_n\to\cuxq$ is injective for each
$n\le\infty$.
In particular, $\pi:\cuoo\to\cuxq$ is injective. Moreover,
this map is surjective because every element $G\in\cuxq$
is the limit of some sequence 
$(G_n)$ of graphs in $\cU$ with $v(G_n)\to\infty$; by \refT{T1},
this sequence converges in $\cuq$ to some element $G'$, and then $\pi(G')=G$.
Since $\cuoo$ is compact, the restriction of $\pi$ to $\cuoo$ is thus
a homeomorphism, and we have the following theorem,
saying that we can identify the set $\cuoo$ of proper graph limits with
$\cuxq$.

\begin{theorem}\label{Tcuoo}
The projection $\pi$ maps  
the set $\cuoo\=\cuq\setminus\cU$ of 
proper graph limits homeomorphically onto
$\cuxq$.
\end{theorem}

\section{Convergence of random graphs}\label{Sconv}

A \emph{random unlabelled graph} is a random element of $\cU$ (with
any distribution; we do not imply any particular model).
We consider convergence of a sequence $(G_n)$ of random unlabelled
graphs in the larger space $\cU$; recall that this is a compact metric
space so we may use the general theory set forth in, for example,
Billingsley \cite{Bill}.

We use the standard notations $\dto$, $\pto$, $\asto$ for convergence
in distribution, probability, and alsmost surely, respectively.
We will only consider the case when $v(G_n)\to\infty$, at least in
probability. (The reader may think of the case when $G_n$ has $n$
vertices, although that is not necessary in general.)

We begin with convergence in distribution.

\begin{theorem}
  \label{T2}
Let $G_n$, $n\ge1$, be random unlabelled graphs and assume that
$v(G_n)\pto\infty$.
The following are equivalent, as \ntoo.
\begin{romenumerate}
  \item\label{T2a}
$G_n\dto \gG$ for some random $\gG\in\cuq$.
  \item\label{T2b}
For every finite family $F_1,\dots,F_m$ of (non-random) graphs, the
random variables $t(F_1,G_n),\dots,t(F_m,G_n)$ converge jointly in
distribution. 
  \item\label{T2c}
For every (non-random) $F\in\cU$, the
random variables $t(F,G_n)$ converge in
distribution.
  \item\label{T2d}
For every (non-random) $F\in\cU$, the
expectations $\E t(F,G_n)$ converge.
\end{romenumerate}
If these properties hold, then the limits in \ref{T2b},
\ref{T2c} and \ref{T2d} are $\bigpar{t(F_i,\gG)}_{i=1}^m$, $t(F,\gG)$ and
$\E t(F,\gG)$, respectively. 
Furthermore, $\gG\in\cuoo$ a.s.

The same results hold if $t$ is replaced by $\tinj$ or $\tind$.
\end{theorem}

\begin{proof}
  \ref{T2a}$\iff$\ref{T2b}.
Since $\cuq$ is a closed subset of $\oiuu$, convergence in
distribution in $\cuq$ is equivalent to convergence of
$\tauu(G_n)=\bigpar{(t(F,G_n))_{F\in\cU},\,v(G_n)\qw}$ in $\oiuu$,
Since we assume $v(G_n)\qw\pto0$, this is equivalent to convergence of
$(t(F,G_n))_{F\in\cU}$ in $\oiu$ \cite[Theorem 4.4]{Bill},
which is equivalent to convergence in distribution of all finite
families $(t(F_i,G_n))_{i=1}^m$.

  \ref{T2b}$\implies$\ref{T2c}. Trivial.

  \ref{T2c}$\implies$\ref{T2d}. Immediate, since $t$ is
  bounded (by 1).

  \ref{T2d}$\implies$\ref{T2b}.
Let $F_1,\dots,F_m$ be fixed graphs and let $\ell_1,\dots,\ell_m$ be
positive integers. Let $F$ be the disjoint union of $\ell_i$ copies of
$F_i$, $i=1,\dots,m$. Then, for every $G\in\cU$, from the definition
of $t$,
\begin{equation*}
  t(F,G)=\prod_{i=1}^mt(F_i,G)^{\ell_i},
\end{equation*}
and hence
\begin{equation}\label{b3}
\E\prod_{i=1}^mt(F_i,G)^{\ell_i}
=\E  t(F,G).
\end{equation}
Consequently, if \ref{T2d} holds, then every joint moment
$\E\prod_{i=1}^mt(F_i,G)^{\ell_i}$ of 
$t(F_1,G_n),\dots,t(F_m,G_n)$ converges. Since $t(F_i,G_n)$ are
bounded (by 1), this implies joint convergence in distribution by the
method of moments.

The identification of the limits is immediate.
Since $v(G_n)\pto\infty$, \ref{T2a} implies that $v(\gG)=\infty$ a.s.,
and thus $\gG\in\cuoo$.

Finally, it follows from \eqref{a4}, \eqref{a4a} and \eqref{a4b} that
we can replace $t$ by $\tinj$ or $\tind$ in \ref{T2b} and \ref{T2d},
and the implications
\ref{T2b}$\implies$\ref{T2c} and \ref{T2c}$\implies$\ref{T2d} are
immediate for $\tinj$ and $\tind$ too.
\end{proof}

Specializing to the case of a non-random limit $G\in\cuoo$, we obtain the
corresponding result for convergence in probability.

\begin{corollary}
  \label{C2}
Let $G_n$, $n\ge1$, be random unlabelled graphs such that
$v(G_n)\pto\infty$, and let $G\in\cuoo$.
The following are equivalent, as \ntoo.
\begin{romenumerate}
  \item\label{C2a}
$G_n\pto G$.
  \item\label{C2c}
$t(F,G_n)\pto t(F,G)$ for every (non-random) $F\in\cU$.
  \item\label{C2d}
$\E t(F,G_n)\to t(F,G)$ for every (non-random) $F\in\cU$.
\end{romenumerate}

The same result holds if $t$ is replaced by $\tinj$ or $\tind$.
\end{corollary}

Note further that under the same assumptions, 
it follows directly from \refT{T1} that
$G_n\asto G$ if and only if
$t(F,G_n)\asto t(F,G)$ for every $F\in\cU$.

We observe another corollary to \refT{T2} (and its proof).

\begin{corollary}
  \label{C2aa}
If $\gG$ is a random element of $\cuoo=\cuq\setminus\cU\cong\cuxq$,
then, for every sequence $F_1,\dots,F_m$ of graphs, possibly with
repetitions,
\begin{equation}\label{b3x}
\E\prod_{i=1}^mt(F_i,\gG)
=\E  t\lrpar{\oplus_{i=1}^mF_i,\gG},
\end{equation}
where $\oplus_{i=1}^mF_i$ denotes the disjoint union of
$F_1,\dots,F_m$. As a consequence, the distribution of $\gG$ is uniquely
determined by the numbers $\E t(F,\gG)$, $F\in\cU$.
Alternatively, 
the distribution of $\gG$ is uniquely
determined by the numbers $\E \tind(F,\gG)$, $F\in\cU$.
\end{corollary}

\begin{proof}
  Since $\cU$ is dense in $\cuq\supseteq\cuoo$, there exists random
  unlabelled graphs $G_n$ such that $G_n\asto \gG$. In particular,
$G_n\dto \gG$ and $v(G_n)\pto\infty$ (in fact, we may assume  $v(G_n)=n$),
so \refT{T2} and its proof apply, and \eqref{b3x} follows from
  \eqref{b3} applied to $G_n$ by letting \ntoo.

For the final statement, note that \eqref{b3x} shows that the
expectations $\E t(F,\gG)$, $F\in\cU$, determine all moments
$\E\prod_{i=1}^m t(F_i,\gG)$, and thus the joint distribution of
$t(F,\gG)$, $F\in\cU$, which is the same as the distribution of
$\tau(\gG)=\bigpar{t(F,\gG)}_{F\in\cU}\in\oiu$, and we have defined
$\cuoo$ such that we identify $\gG$ and $\tau(\gG)$.
Finally, the numbers $\E \tind(F,\gG)$, $F\in\cU$,
determine all $\E t(F,\gG)$ by 
\eqref{a4a}, recalling that $\tinj(F,\gG)=t(F,\gG)$ by \eqref{a4x}.
 \end{proof}

\begin{remark}
  The numbers $\E t(F,\gG)$ for a random $\gG\in\cuoo$ thus play a role
  similar to the one played by moments for a random variable.
(And the relation between $\E t(F,\gG)$ and $\E\tind(F,\gG)$ has some
  resemblance to the relation between moments and cumulants.)
\end{remark}

\section{Convergence to infinite graphs}\label{Sinfinite}

We will in this section consider also labelled \emph{infinite} graphs
with the vertex set $\bbN=\set{1,2,\dots}$.
Let $\cloo$ denote the set of all such graphs. These graphs are
determined by their edge sets, so $\cloo$ can be identified with the power
set $\cP(E(\Koo))$
of all subsets of the edge set $E(\Koo)$ of the complete infinite
graph $\Koo$, and thus with the infinite product set $\set{0,1}^{E(\Koo)}$.
We give this space, and thus $\cloo$, the product topology. Hence,
$\cloo$ is a compact metric space.

It is sometimes convenient to regard $\cL_n$ for a finite $n$ as a
subset of $\cloo$: we can identify graphs in $\cL_n$ and $\cloo$ with
the same edge set. In other words, if $G\in\cL_n$ is a graph with
vertex set $[n]$, we add an infinite number of isolated vertices
$n+1,n+2,\dots$ to obtain a graph in $\cloo$.

Conversely, if $H\in\cloo$ is an infinite graph, we let $H\rest
n\in\cL_n$ be the induced subgraph of $H$ with vertex set $[n]$.

If $G$ is a (finite) graph, let $\hG$ be the random labelled graph
obtained by a random labelling of the vertices of $G$ by the numbers
$1,\dots,v(G)$.
(If $G$ is labelled, we thus ignore the labels and randomly relabel.)
Thus $\hG$ is a random finite graph with the same number of vertices
as $G$, but  as just said, we can (and will) also regard $\hG$ as a
random graph in $\cloo$.

We use the same notation $\hG$ also for a random (finite) graph $G$
given a random labelling.

\begin{theorem}\label{TC1}
  Let $(G_n)$ be a sequence of random graphs in $\cU$ and assume that
  $v(G_n)\pto\infty$. Then the following are equivalent.
  \begin{romenumerate}
\item
$G_n\dto \gG$ in $\cuq$ for some random $\gG\in\cuq$.	
\item
$\hgn\dto H$ in $\cloo$ for some random $H\in\cloo$.	
  \end{romenumerate}
If these hold, then
$\P\xpar{H\rest k =F} = \E \tind(F,\gG)$ for every $F\in\cL_k$.
Furthermore, $\gG\in\cuoo$ a.s.
\end{theorem}

\begin{proof}
Let $G$ be a labelled graph and consider the graph $\hG\rest k$,
assuming $k\le v(G)$. This random graph equals $G\xk'=G(\vvki)$, where
$\vvki$ are $k$ vertices sampled at random without replacement as in
\refS{Sdef}. 
Hence, by \eqref{tind}, for every $F\in\cL_k$,
\begin{equation*}
  \P\xpar{\hG\rest k =F} = \tind(F,G),
\qquad \text{if } k\le v(G).
\end{equation*}
Applied to the random graph $G_n$, this yields
\begin{equation}
  \label{c3}
\E\tind(F,G_n)
\le
\P\xpar{\hgn\rest k =F} \le \E\tind(F,G_n) + P\bigpar{v(G_n)<k}.
\end{equation}
By assumption, $P\lrpar{v(G_n)<k}\to0$ as \ntoo, and it follows from
\eqref{c3} and \refT{T2} that $G_n\dto \gG$ in $\cuq$ if and only if
\begin{equation}\label{c3a}
  \P\xpar{\hgn\rest k =F} \to \E\tind(F,\gG)
\end{equation}
for every $k\ge1$ and every $F\in\cL_k$.

Since $\cL_k$ is a finite set, \eqref{c3a} says that, for every $k$,
$\hgn\rest k\dto H_k$ for some random graph $H_k\in\cL_k$ with
$\P(H_k=F)=\E\tind(F,\gG)$ for $F\in\cL_k$. Since $\cloo$ has the
product topology, this implies $\hgn\dto H$ in $\cloo$ for some random
$H\in\cloo$ with $H\rest k\eqd H_k$.

Conversely, if $\hgn\dto H$ in $\cloo$, then $\hgn\rest k\dto H\rest
k$ so the argument above shows that
\begin{equation*}
  \E\tind(F,G_n)=\P\xpar{\hgn\rest k=F}+o(1)\to\P\xpar{H\rest k=F}
\end{equation*}
as \ntoo, for every $F\in\cL_k$, and \refT{T2} yields the existence of
some random $\gG\in\cuoo\subset\cuq$ with $G_n\dto \gG$
and
$\E\tind(F,\gG)=\P\xpar{H\rest k=F}$. 
\end{proof}

\section{Exchangeable random graphs}\label{Sexch}

\begin{definition}
  A random infinite graph $H\in\cloo$ is \emph{\exch} if its
  distribution is invariant under every permutation of the vertices.
(It is well-known that it is equivalent to consider only finite 
  permutations, \ie, permutations $\gs$ of $\bbN$ that satisfy
  $\gs(i)=i$ for all sufficiently large $i$, so $\gs$ may be regarded
  as a permutation in $\sn$ for some $n$.)

Equivalently, if $X_{ij}\=\ett{ij\in H}$ is the indicator of there 
being an edge $ij$ in $H$, then the array \set{X_{ij}}, $1\le
i,j\le\infty$, is (jointly) \exch{} as defined in \refS{S:intro}.
\end{definition}

\begin{lemma}
  \label{LE}
Let $H$ be a random infinite graph in $\cloo$. Then the following are
equivalent. 
\begin{romenumerate}
  \item
$H$ is \exch.
\item
$H\rest k$ has a distribution invariant under all permutations of
$[k]$, for every $k\ge1$.
\item
$\P\bigpar{H\rest k=F}$ depends only on the isomorphism type of $F$,
and can thus be seen as a function of $F$ as an unlabelled graph in
$\cU_k$, for every $k\ge1$.
\end{romenumerate}
\end{lemma}
\begin{proof}
  (i)$\implies$(ii). Immediate.

  (ii)$\implies$(i). If $\gs$ is a finite permutation of $\bbN$, then
$\gs$ restricts to a permutation of $[k]$ for every large $k$, and it
  follows that if $\Hgs$ is $H$ with the vertices permuted by $\gs$,
  then, for all large $k$
$\Hgs\rest k={H\rest k}\permgs\eqd H\rest k$, which implies 
$\Hgs\eqd H$.

(ii)$\iff$(iii). Trivial.
\end{proof}

\begin{theorem}
  \label{TC1E}
  The limit $H$ is \refT{TC1} is \exch.
\end{theorem}
\begin{proof}
$H$ satisfies \refL{LE}(iii).  
\end{proof}

Moreover, \refT{TC1} implies the following
connection with random elements of $\cuoo$.

\begin{theorem}
  \label{TE}
There is a one-to-one correspondence between distributions of random
elements $\gG\in\cuoo$ (or $\cuxq$) and 
distributions of
\exch{} random infinite graphs
$H\in\cloo$ given by
\begin{align}
  \label{e2a}
\E\tind(F,\gG)&=\P(H\rest k=F)
\intertext{for every $k\ge1$ and every $F\in\cL_k$, or, equivalently,}
  \label{e2b}  
\E t(F,\gG)&=\P(H\supset F)
\end{align}
for every $F\in\cL$. 
Furthermore, $H\rest n\dto \gG$ in $\cuq$ as \ntoo.
\end{theorem}

\begin{proof}
Note first that \eqref{e2a} and \eqref{e2b} are equivalent by
\eqref{a4a} and \eqref{a4b}, since $t(F,\gG)=\tinj(F,\gG)$ by \eqref{a4x},
and $H\supset F$ if and only if $H\rest k\supseteq F$ when $F\in\cL_k$.

Suppose that $\gG$ is a random element of
  $\cuoo\subset\cuq$. Since $\cU$ is dense in $\cuq$, there exist (as
  in the proof of \refC{C2aa})
  random unlabelled graphs $G_n$ such that $G_n\asto \gG$ in $\cuq$ and
  thus $v(G_n)\asto\infty$ and $G_n\dto \gG$.
Hence, Theorems \refand{TC1}{TC1E} 
show that $\hgn\dto H$ for some random \exch{} infinite graph $H$
  satisfying \eqref{e2a}.
Furthermore, \eqref{e2a} determines the distribution of $H\rest k$ for
  every $k$, and thus the distribution of $k$.

Conversely, if $H$ is an \exch{} random infinite graph, let
$G_n=H\rest n$. By \refL{LE}(ii), the distribution of each $G_n$ is
invariant under permutations of the vertices, so if $\hgn$ is $G_n$
with a random (re)labelling, we have $\hgn\eqd G_n$.
Since $G_n\dto H$ in $\cloo$ (because $\cloo$ has a product topology),
we thus have $\hgn\dto H$ in $\cloo$, so \refT{TC1} applies and shows
the existence of a random $\gG\in\cuoo$ such that 
$G_n\dto \gG$ and 
\eqref{e2a} holds. 
Finally \eqref{e2a} determines the distribution of $\gG$ 
by \refC{C2aa}.
\end{proof}

\begin{remark}\label{Ras}
  Moreover, $H\rest n$ converges \as{} to some random variable $\gG\in\cuoo$, 
because $\tind(F,H\rest n)$, $n\ge v(F)$, 
is a reverse martingale for every $F\in \gG$. 
Alternatively, this follows by concentration estimates from the
representation in \refS{Srep}, see
\Lovaszetal{} \cite[Theorem 2.5]{LSz}.
\end{remark}

\begin{corollary}\label{CE}
There is a one-to-one correspondence between 
elements $\gG$ of $\cuoo\cong\cuxq$ and extreme points of the set of
distributions of \exch{} random infinite graphs
$H\in\cloo$.
This correspondence is given by 
\begin{equation}\label{ce}
t(F,\gG)=\P(H\supset F)
\end{equation}
for every $F\in\cL$. Furthermore, $H\rest n\asto \gG$ in $\cuq$ as \ntoo.
\end{corollary}
\begin{proof}
  The extreme points of the set of distributions on $\cuoo$ are the
  point masses, which are in one-to-one correspondence with the
  elements of $\cuoo$.
\end{proof}

We can characterize these extreme point distributions of \exch{}
random infinite graphs as follows.

\begin{theorem}
  \label{TE2}
Let $H$ be an \exch{} random infinite graph.
Then the following are equivalent.
\begin{romenumerate}
  \item\label{te2a}
The distribution of $H$ is an extreme point in the set of \exch{}
distributions in $\cloo$.
  \item\label{te2b'}
If $F_1$ and $F_2$ are two (finite) graphs with disjoint vertex sets
$V(F_1)$, $V(F_2)\subset\bbN$, then
\begin{equation*}
  \P(H\supset F_1\cup F_2)=\P(H\supset F_1)\P(H\supset F_2).
\end{equation*}
  \item\label{te2b}
The restrictions $H\rest k$ and $H\restx{[k+1,\infty)}$ are
  independent for every $k$.
  \item\label{te2c}
Let $\cF_n$ be the $\gs$-field generated by
$H\restx{[n,\infty)}$. Then the tail $\gs$-field
  $\bigcap_{n=1}^\infty \cF_n$ is trivial, \ie, contains only events with
  probability $0$ or $1$.
\end{romenumerate}
\end{theorem}
\begin{proof}
  (i)$\implies$(ii).
By \refC{CE}, $H$ corresponds to some (non-random) $\gG\in\cuoo$ such
that
\begin{equation}
  \label{e6}
\P(H\supset F)=t(F,\gG)
\end{equation}
for every $F\in\cL$. We have defined $\cL$ such that a graph $F\in\cL$
is labelled by $1,\dots,v(F)$, but both sides of \eqref{e6} are
invariant under relabelling of $F$ by arbitrary positive integers; the
\lhs{} because $H$ is \exch{} and the \rhs{} because $t(F,\gG)$ only
depends on $F$ as an unlabelled graph. Hence \eqref{e6} holds for
every finite graph $F$ with $V(F)\subset\bbN$.

Furthermore, since $\gG$ is non-random, \refC{C2aa} yields 
$t(F_1\cup F_2,\gG)=t(F_1,\gG)t(F_2,\gG)$. Hence,
\begin{equation*}
  \begin{split}
\P(H\supset F_1\cup F_2)=
t(F_1\cup F_2,\gG)
=t(F_1,\gG)t(F_2,\gG)
=\P(H\supset F_1)\P(H\supset F_2).
  \end{split}
\end{equation*}

  (ii)$\implies$(iii).
By inclusion--exclusion, as for \eqref{tind}, 
(ii) implies that 
if $1\le k < l<\infty$, then 
for any graphs $F_1$ and $F_2$ with 
$V(F_1)=\set{1,\dots,k}$ and $V(F_2)=\set{k+1,\dots,k+l}$,
the events $H\rest k=F_1$ and $H\restx{\set{k+1,\dots,l}} =F_2$ 
are independent.
Hence $H\rest k$ and $H\restx{\set{k,\dots,l}}$ are independent for
every $l>k$, and the result follows.

  (iii)$\implies$(iv).
Suppose $A$ is an event in the tail $\gs$-field
  $\bigcap_{n=1}^\infty \cF_n$.
Let $\cF^*_n$ be the $\gs$-field generated by $H\rest n$. By
\ref{te2b}, $A$ is independent of $\cF^*_n$ for every $n$, and thus of
the $\gs$-field $\cF$ generated by $\bigcup\cF^*_n$, which equals the
$\gs$-field $\cF_1$ generated by $H$. However, $A\in\cF_1$, so $A$ is
independent of itself and thus $\P(A)=0$ or 1.

  (iv)$\implies$(i).
Let $F\in\cL_k$ for some $k$ and let $F_n$ be $F$ with all vertices
shifted by $n$. Consider the two indicators $I=\ett{H\supseteq F}$ and
$I_n=\ett{H\supseteq F_n}$. Since $I_n$ is $\cF_n$-measurable,
\begin{equation}\label{e70}
  \P(H\supset F\cup F_n)
=\E(II_n)
=\E\bigpar{\E(I\mid\cF_n)I_n}.
\end{equation}
Moreover, $\E(I\mid\cF_n)$, $n=1,2,\dots$, is a reverse martingale,
and thus a.s.  
\begin{equation*}
\E(I\mid\cF_n)\to\E\Bigpar{I\mid\bigcap_{n=1}^\infty \cF_n} 
=\E I,
\end{equation*}
using (iv). 
Hence, $\bigpar{\E(I\mid\cF_n)-\E I}I_n\to0$ a.s., and by dominated
convergence 
\begin{equation*}
\E\Bigpar{ \bigpar{\E(I\mid\cF_n)-\E I}I_n}\to0.
\end{equation*}
Consequently, \eqref{e70} yields
\begin{equation*}
  \P(H\supset F\cup F_n) = \E I \E I_n + o(1)
=
\P(H \supset F) \P(H\supset F_n) + o(1).
\end{equation*}
Moreover, since $H$ is \exch, $\P(H\supset F\cup F_n)$ (for $n\ge v(F)$)
and $\P(H\supset F_n)$ do not depend on $n$, and we
obtain as \ntoo
\begin{equation}\label{e7}
  \P(H\supset F\cup F_k)=\P(H \supset F)^2.
\end{equation}
Let $\gG$ be a random element of $\cuoo$ corresponding to $H$ as in
\refT{TE}.
By \eqref{e2b} and \eqref{b3x}, \eqref{e7} can be written
\begin{equation*}
  \E t(F,\gG)^2=\bigpar{\E t(F,\gG)}^2.
\end{equation*}
Hence the random variable $t(F,\gG)$ has variance 0 so it is \as{}
constant. Since this holds for every $F\in\cL$, it follows that $\gG$ is
\as{} constant, \ie, we can take $\gG$ non-random, and \ref{te2a}
follows
by \refC{CE}.
\end{proof}

\section{Representations of graph limits and \exch{} graphs}\label{Srep}

As said in the introduction, 
the \exch{} infinite random graphs were characterized by Aldous
\cite{Aldous}
and Hoover \cite{Hoover},
see also Kallenberg \cite{Kallenberg:exch}, and the graph limits in
$\cuoo\cong\cuxq$ were characterized in a very similar way 
by \Lovaszetal{} \cite{LSz}. We can now make the connection between
these two characterizations explicit.

Let $\cw$ be the set of all measurable functions $W:\oi^2\to\oi$ and
let $\cws$ be the subset of symmetric functions.
For every $W\in\cws$, we define an infinite random graph
$\gwoo\in\cloo$ as follows: we
first choose a sequence $X_1,X_2,\dots$ of \iid{} random variables
uniformly distributed on $\oi$, and then, given this sequence,
for each pair $(i,j)$ with $i<j$ we draw an edge 
$ij$ with probability $W(X_i,X_j)$, independently for
all pairs $(i,j)$ with $i<j$ (conditionally given \set{X_i}).
Further, let $\gwn$ be the restriction $\gwoo\rest n$, which is
obtained by the same construction with a finite sequence $X_1,\dots,X_n$.

It is evident that $\gwoo$ is an \exch{} infinite random graph, and
the result by Aldous and Hoover 
is that every \exch{} infinite
random graph is obtained as a mixture of such $\gwoo$; in other words
as $\gwoo$ with a random $W$.

Considering again a deterministic $W\in\cws$, it is evident that
\refT{TE2}(ii) holds, and thus \refT{TE2} and \refC{CE} show that
$\gwoo$ 
corresponds to an element $\uw\in\cuoo$.
Moreover, by \refT{TE} and \refR{Ras}, $\gwn\to \uw$ \as{} as \ntoo, 
and \eqref{ce}
shows that if $F\in\cL_k$, then
\begin{equation}\label{tuw}
  t(F,\uw)=
\P\bigpar{F\subseteq \gwx k}
=
\int_{\oi^k} \prod_{ij\in\E(F)} W(x_i,x_j) \dd x_1\dots \dd x_k.
\end{equation}
The main result of \Lovaszetal{} \cite{LSz} is that every element of
$\cuoo\cong \cuxq$ can be obtained as $\uw$ satisfying \eqref{tuw} for
some $W\in\cws$.

It is now  clear that the representation theorems of Aldous--Hoover
\cite{Aldous,Hoover} and \Lovaszetal{} \cite{LSz} are connected by 
\refT{TE} and \refC{CE} above, and that one characterization easily
follows from the other.

\begin{remark}
  The representations by $W$ are far from unique, see \refS{Sunique}.
\Borgsetal{} \cite{BCL1} call an element $W\in\cws$ a \emph{graphon}.
They further define a pseudometric (called the
\emph{cut-distance}) on $\cws$ and show that if we consider the quotient
space $\cwsq$ obtained by identifying elements with cut-distance 0,
we obtain a compact metric space, and the mapping $W\mapsto\uw$ yields
a bijection $\cwsq\to\cuxq\cong\cuoo$, which furthermore is a homeomorphism.
\end{remark}

\begin{remark}
As remarked in \Lovaszetal{} \cite{LSz}, we can more generally 
consider a symmetric measurable function $W:\cS^2\to\oi$ for any
probability space $(\cS,\mu)$, and define $\gwoo$ 
as above with $X_i$ \iid{} random variables in $\cS$ with distribution
$\mu$. This does not give any new limit objects $\gwoo$ or $\uw$,
since we just said that every limit object is obtained from some $W\in\cws$,
but they can sometimes give useful representations.

An interesting case is when $W$ is the adjacency matrix of a (finite)
graph $G$, with $\cS=V(G)$ and $\mu$ the uniform measure on $\cS$; we
thus let $X_i$ be \iid{} random vertices of $G$ and $\gwn$ equals the
random graph $G[n]$ defined in \refS{Sdef}. It follows from \eqref{tuw} and
\eqref{t} that $t(F,\uw)=t(F,G)$ for every $F\in\cU$, and thus $\uw=G$
as elements of $\cuxq$. In other words, $\uw\in\cuoo=\pi(G)$, the
``ghost'' of $G$ in $\cuoo\cong\cuxq$.
\end{remark}

\begin{remark}
For the asymptotic behavior of $G(n,W)$ in another, sparse, case,
with $W$ depending on $n$, see \cite{SJ178}.
\end{remark}

\section{Non-uniqueness}\label{Sunique}

The functions $W$ on $\oi^2$ used to represent graph limits or \exch{}
arrays are far from unique.
(For a special case when there is a natural canonical choice, which much
simplifies and helps applications, see \cite{threshold}.)
For example, it is obvious that if $\gf:\oi\to\oi$ is any measure
preserving map, then $W$ and $\W\gf$, defined by
$\W\gf(x,y)\=W\bigpar{\gf(x),\gf(y)}$, define the same graph limit 
and the same (in distribution) \exch{} array.

Although in principle, this is the only source on non-uniqueness, the
details are more complicated, mainly because the measure preserving
map $\gf$ does not have to be a bijection, and thus the relation $W'=\W\gf$
is not symmetric: it can hold without there being a measure preserving
map $\gf'$ such that $W=W'\circ \gf'$. (For a 1-dimensional example,
consider $f(x)=x$ and  $f'(x)=\gf(x)=2x \bmod1$; for a 2-dimensional
example, let $W(x,y)=f(x)f(y)$ and $W'(x,y)=f'(x)f'(y)$.)

For \exch{} arrays, the equivalence problem was solved by
\citet{Hoover}, who gave a criterion which in our case reduces to \ref{TUpsi}
below; this criterion involves an auxiliary variable, and can
be interpreted as saying $W=W'\circ\gf'$ for a random $\gf'$.
This work was continued by Kallenberg, see \cite{Kallenberg:exch}, who
gave a probabilistic proof and added criterion \ref{TUphi}. For graph
limits, 
\citet{BCL1} gave the criterion \ref{TUcut} in terms of the cut-distance,
and \citet{BR} found the criterion \ref{TUphi} in this context.
Further,
\citet{BCL1} announced the related criterion that there exists a 
measurable function $U:\oi^2\to\oi$ and two measure preserving maps
$\gf,\gf':\oi\to\oi$ such that $W=U\circ\gf$ and $W'=U\circ\gf'$ a.e.;
the proof of this will appear in \cite{BCL:unique}.

As in \refS{Srep}, these two lines of work are connected by the
results in \refS{Sexch}, and we can combine the previous results as
follows.

\begin{theorem}\label{TU}
  Let $W,W'\in\cws$. Then the following are equivalent.
  \begin{romenumerate}
\item\label{TUgg}
$\gG_W=\gG_{W'}$ for the graph limits $\gG_W,\gG_{W'}\in\cuoo$.	
\item\label{TUt}
$t(F,\gG_W)=t(F,\gG_{W'})$ for every graph $F$.
\item\label{TUgoo}
The \exch{} random infinite graphs $G(\infty,W)$ and $G(\infty,W')$
have the same distribution.
\item\label{TUgn}
The random graphs $G(n,W)$ and $G(n,W')$
have the same distribution for every finite $n$.
\item\label{TUphi}
There exist measure preserving maps $\gf,\gf':\oi\to\oi$ such that
$\W\gf=\W{\gf'}$ \aex{} on $\oi^2$, \ie, 
$W\bigpar{\gf(x),\gf(y)}=W'\bigpar{\gf'(x),\gf'(y)}$
a.e.
\item\label{TUpsi}
There exists a measure preserving map $\psi:\oi^2\to\oi$ such that
$W(x_1,x_2)=W'\bigpar{\psi(x_1,y_1),\psi(x_2,y_2)}$ \aex{} on $\oi^4$.
\item\label{TUcut}
$\dcut(W,W')=0$, where $\dcut$ is the cut-distance defined in \cite{BCL1}.
  \end{romenumerate}
\end{theorem}

\begin{proof}

\ref{TUgg}$\iff$\ref{TUt}. By our definition of $\cuoo\subset\cuq$.

\ref{TUgg}$\iff$\ref{TUgoo}. By \refC{CE}.

\ref{TUgoo}$\iff$\ref{TUgn}. Obvious.

\ref{TUphi}$\implies$\ref{TUgoo}. If $X_1,X_2,\dots$ are \iid{} random
variables 
uniformly distributed on $\oi$,
then so are $\gf(X_1),\gf(X_2),\dots$, and thus
$G(\infty,W)\eqd G(\infty,\W\gf) = G(\infty,W'\circ\gf')
\eqd G(\infty,W')$.

\ref{TUgoo}$\implies$\ref{TUphi}. 
The general form of the representation theorem as stated in
\cite[Theorem 7.15, see also p.\ 304]{Kallenberg:exch} is (in our
two-dimensional case)
$X_{ij}=f(\xio,\xii,\xij,\xiij)$ for a function $f:\oi^4\to\oi$,
symmetric in the two middle variables, and 
independent random variables 
$\xio$, $\xii$ ($1\le i$) and $\xiij$ ($1\le i<j$), all  
uniformly distributed on \oi, and where we further let 
$\xi_{ji}=\xiij$ for $j>i$.
We can write the construction of $G(\infty,W)$ in this form with
\begin{equation}
  \label{fW}
f(\xio,\xii,\xij,\xiij)=\ett{\xiij\le W(\xii,\xij)}.
\end{equation}
Note that this $f$ does not depend on $\xio$. (In general, $\xio$ is
needed for the case of a random $W$, which can be written as a
deterministic function of $\xio$, but this is not needed in the
present theorem.) 

Suppose that $G(\infty,W)\eqd G(\infty,W')$, let $f$ be given by $W$
by \eqref{fW}, and let similarly $f'$ be given by $W'$; for notational
convenience we write 
$W_1\=W$, $W_2\=W'$, $f_1\=f$ and $f_2\=f'$.
The equivalence theorem \cite[Theorem 7.28]{Kallenberg:exch}
takes the form, using \eqref{fW}, that there exist measurable
functions
$\gko:\oi\to\oi$, $\gki:\oi^2\to\oi$ and $\gkii:\oi^4\to\oi$, for $k=1,2$, 
that are measure preserving in the last coordinate for any fixed
values of the other coordinates, and such that the two functions
(for $k=1$ and $k=2$)
\begin{multline*}
f_k\bigpar{\gko(\xio),\gki(\xio,\xi_1),\gki(\xio,\xi_2),
 \gkii(\xio,\xi_1,\xi_2,\xi_{12})}
\\
=
\bigett{
W_k\bigpar{\gki(\xio,\xi_1),\gki(\xio,\xi_2)} \ge
 \gkii(\xio,\xi_1,\xi_2,\xi_{12})}
\end{multline*}
are \as{} equal.
Conditioned on $\xio,\xi_1$ and $\xi_2$, the random variable 
$ \gkii(\xio,\xi_1,\xi_2,\xi_{12})$ is uniformly distributed on $\oi$,
and it follows (e.g., by taking the conditional expectation) that \as
\begin{equation*}
  W_1\bigpar{\gii(\xio,\xi_1),\gii(\xio,\xi_2)} 
=
W_2\bigpar{\gji(\xio,\xi_1),\gji(\xio,\xi_2)}.
\end{equation*}
For \aex{} value $x_0$ of $\xio$, this thus holds for \aex{} values of $\xi_1$
and $\xi_2$, and we may choose $\gf(x)=\gii(x_0,x)$ and $\gf'(x)\=\gji(x_0,x)$
for some such $x_0$.

\ref{TUgoo}$\iff$\ref{TUpsi}. 
Similar, using \cite[Theorem 7.28(iii)]{Kallenberg:exch}.

\ref{TUt}$\iff$\ref{TUcut}. See \cite{BCL1}.
\end{proof}

\section{Bipartite graphs}\label{Sbip}

The definitions and results above have analogues for bipartite
graphs, which we give in this section, leaving some details to the reader.
The proofs are straightforward analogues of the ones given above and
are omitted.
Applications of the results of
this section to random difference graphs are in 
\cite{threshold}.

A \emph{bipartite graph} will be a graph with an explicit bipartition; in
other words, a bipartite graph $G$ consists of two vertex sets
$V_1(G)$ and $V_2(G)$ and an edge set $E(G)\subseteq V_1(G)\times
V_2(G)$; we let $v_1(G)\=|V_1(G)|$ and $v_2(G)\=|V_2(G)|$ be the
numbers of vertices in the two sets.
Again we consider both the labelled and unlabelled cases; in the
labelled case we assume the labels of the vertices in $V_j(G)$ are
$1,\dots,v_j(G)$ for $j=1,2$.
Let $\BLnn$ be the set of the $2^{n_1n_2}$ labelled bipartite graphs
with vertex sets $[n_1]$ and $[n_2]$, and let $\BUnn$ be the quotient
set $\BLnn/\cong$ of unlabelled bipartite graphs with $n_1$ and $n_2$
vertices in the two parts; further, let
$\BL\=\bigcup_{n_1,n_2\ge1}\BLnn$ and $\BU\=\bigcup_{n_1,n_2\ge1}\BUnn$.

We let $G[k_1,k_2]$ be the random graph in $\BLxx{k_1}{k_2}$ obtained
by sampling $k_j$ vertices from $V_j(G)$ ($j=1,2$), uniformly with
replacement, and let, provided $k_j\le v_j(G)$, 
$G[k_1,k_2]'$ be the corresponding random graph obtained by sampling
without replacement. We then define
$t(F,G)$, $\tinj(F,G)$ and $\tind(F,G)$ 
for (unlabelled) bipartite graphs $F$ and $G$
in analogy with \eqref{t}--\eqref{tind}.
Then \eqref{a4}--\eqref{a4b} still hold, \emph{mutatis mutandis}; for
example,
\begin{equation}
  \label{a4xx}
\bigabs{t(F,G)-\tinj(F,G)}
\le\frac{v_1(F)^2}{2v_1(G)}+\frac{v_2(F)^2}{2v_2(G)}.
\end{equation}

In analogy with \eqref{tau}, we now define $\tau:\BU\to\oib$ by
\begin{equation}\label{taub}
  \tau(G)\=(t(F,G))_{F\in\cB}\in\oib. 
\end{equation}

We define $\cbx\=\tau(\cB)\subseteq\oib$ to be the image of $\cB$
under this mapping $\tau$, and let $\cbxq$ be the closure of $\cbx$ in
$\oib$; this is
a compact metric space.

Again, $\tau$ is not injective; we may
consider a graph $G$ as an element of $\cbx$ by identifying $G$ and
$\tau(G)$, but this implies identification of some
graphs of different orders and 
we prefer to avoid it.
We let
$\cbb$ be the union of $\cB$ and some two-point set \set{*_1,*_2} and
consider the mapping $\tauu:\cB\to\oibb=\oib\times\oi\times\oi$ defined by
\begin{equation}\label{tab+}
  \tauu(G)=\bigpar{\tau(G),\,v_1(G)\qw,\,v_2(G)\qw}.
\end{equation}
Then $\tauu$ is injective
and we can identify $\cB$ with its image $\tauu(\cB)\subseteq\oibb$
and define $\cbq\subseteq\oibb$ as its closure; this is a compact
metric space.

The functions $t(F,\cdot)$,
$\tinj(F,\cdot)$, $\tind(F,\cdot)$ and 
$v_j(\cdot)\qw$, for $F\in\cB$ and $j=1,2$, 
have unique continuous
extensions to $\cbq$.

We let $\cboo\=\set{G\in\cbq:v_1(G)=v_2(G)=\infty}$; 
this is the set of all limit objects of sequences
$(G_n)$ in $\cB$ with $v_1(G_n),v_2(G_n)\to\infty$.
By \eqref{a4xx}, $\tinj(F,G)=t(F,G)$ for every $G\in\cboo$
and every $F\in\cB$. The projection $\pi:\cbq\to\cbxq$ restricts to a
homeomorphism $\cboo\cong\cbxq$.

\begin{remark}
Note that in the bipartite case there are other limit objects too in
$\cbq$; in fact,  $\cbq$ can be partitioned into $\cB$, $\cboo$, and the sets
$\BUxx{n}{\infty}$, $\BUxx{\infty}n$, for $n=1,2,\dots$, where, for
example, 
$\BUxx{n_1}{\infty}$ is the set of limits of sequences $(G_n)$ of
bipartite graphs such that $v_2(G_n)\to\infty$ but $v_1(G_n)=n_1$ is constant.
We will not consider such degenerate limits further here, but we
remark that in the simplest case $n_1=1$, a bipartite graph in $\BLxx{1}{n_2}$
can be identified with a subset of $[n_2]$, and an unlabelled graph in
$\BUxx{1}{n_2}$ thus with a number in $m\in\set{0,\dots,n_2}$, the number of
edges in the graph, and it is easily seen that a sequence of such
unlabelled graphs with $n_2\to\infty$ converges in $\cbq$ if and only
if the proportion $m/n_2$ converges; hence we can identify 
$\BUxx{1}{\infty}$ with the interval \oi.  
\end{remark}

We have the following basic result, \cf{}
\refT{T1}.

\begin{theorem}\label{T1B}
 Let $(G_n)$ be a sequence of bipartite graphs with
 $v_1(G_n)$, $v_2(G_n)\to\infty$.
Then the following are equivalent.
\begin{romenumerate}
  \item\label{T1Bt}
$t(F,G_n)$ converges for every $F\in\cB$.
  \item\label{T1Btinj}
$\tinj(F,G_n)$ converges for every $F\in\cB$.
  \item\label{T1Btind}
$\tind(F,G_n)$ converges for every $F\in\cB$.
  \item
$G_n$ converges in $\cbq$.
\end{romenumerate}
In this case, the limit $G$ of $G_n$ belongs to $\cboo$ and the limits
in \ref{T1Bt}, \ref{T1Btind} and \ref{T1Btind} are $t(F,G)$,
$\tinj(F,G)$ and $\tind(F,G)$.
\end{theorem}

For convergence of random unlabelled bipartite graphs, the results in
\refS{Sconv} hold with trivial changes.

\begin{theorem}
  \label{T2B}
Let $G_n$, $n\ge1$, be random unlabelled bipartite graphs and assume that
$v_1(G_n),v_2(G_n)\pto\infty$.
The following are equivalent, as \ntoo.
\begin{romenumerate}
  \item\label{T2Ba}
$G_n\dto \gG$ for some random $\gG\in\cbq$.
  \item\label{T2Bb}
For every finite family $F_1,\dots,F_m$ of (non-random) bipartite graphs, the
random variables $t(F_1,G_n),\dots,t(F_m,G_n)$ converge jointly in
distribution. 
  \item\label{T2Bc}
For every (non-random) $F\in\cB$, the
random variables $t(F,G_n)$ converge in
distribution.
  \item\label{T2Bd}
For every (non-random) $F\in\cB$, the
expectations $\E t(F,G_n)$ converge.
\end{romenumerate}
If these properties hold, then the limits in \ref{T2Bb},
\ref{T2Bc} and \ref{T2Bd} are $\bigpar{t(F_i,\gG)}_{i=1}^m$, $t(F,\gG)$ and
$\E t(F,\gG)$, respectively. 
Furthermore, $\gG\in\cboo$ a.s.

The same results hold if $t$ is replaced by $\tinj$ or $\tind$.
\end{theorem}

\begin{corollary}
  \label{C2B}
Let $G_n$, $n\ge1$, be random unlabelled bipartite graphs such that
$v_1(G_n),v_2(G_n)\pto\infty$, and let $G\in\cboo$.
The following are equivalent, as \ntoo.
\begin{romenumerate}
  \item\label{C2Ba}
$G_n\pto G$.
  \item\label{C2Bc}
$t(F,G_n)\pto t(F,G)$ for every (non-random) $F\in\cB$.
  \item\label{C2Bd}
$\E t(F,G_n)\to t(F,G)$ for every (non-random) $F\in\cB$.
\end{romenumerate}

The same result holds if $t$ is replaced by $\tinj$ or $\tind$.
\end{corollary}

As above, the distribution of $\gG$ is uniquely
determined by the numbers $\E t(F,\gG)$, $F\in\cB$.

Let $\cbloo$ denote the set of all 
labelled {infinite} bipartite graphs
with the vertex sets $V_1(G)=V_2(G)=\bbN$.
$\cbloo$ is a compact metric space with the natural product topology.

If $G$ is a bipartite graph, let $\hG$ be the random labelled
bipartite graph 
obtained by random labellings of the vertices in $V_j(G)$ by the numbers
$1,\dots,v_j(G)$, for $j=1,2$.
This is a random finite bipartite graph, but
we  can also regard it as a random
element of $\cbloo$ by adding isolated vertices.

\begin{definition}
  A random infinite bipartite graph $H\in\cbloo$ is \emph{\exch} if its
  distribution is invariant under every pair of finite permutations of
  $V_1(H)$ and $V_2(H)$.
\end{definition}

\begin{theorem}\label{TC1B}
  Let $(G_n)$ be a sequence of random graphs in $\cB$ and assume that
  $v_1(G_n),v_2(G_n)\pto\infty$. Then the following are equivalent.
  \begin{romenumerate}
\item
$G_n\dto \gG$ in $\cbq$ for some random $\gG\in\cbq$.	
\item
$\hgn\dto H$ in $\cbloo$ for some random $H\in\cbloo$.	
  \end{romenumerate}
If these hold, then
$\P\xpar{H\restkk =F} = \E \tind(F,\gG)$ for every $F\in\BLkk$.
Furthermore, $\gG\in\cboo$ a.s., and
$H$ is \exch.
\end{theorem}

\begin{theorem}
  \label{TEB}
There is a one-to-one correspondence between distributions of random
elements $\gG\in\cboo$ (or $\cbxq$) and 
distributions of \exch{} random infinite graphs
$H\in\cbloo$ given by
\begin{align}
\E\tind(F,\gG)&=\P\xpar{H\restkk =F}
\intertext{for every $k_1,k_2\ge1$ and every $F\in\BLkk$, or, equivalently,}
\E t(F,\gG)&=\P(H\supset F)
\end{align}
for every $F\in\BL$.
Furthermore, $H\restnn\dto \gG$ in $\cbq$ as $n_1,n_2\to\infty$.
\end{theorem}

\begin{corollary}\label{CEB}
There is a one-to-one correspondence between 
elements $\gG$ of $\cboo\cong\cbxq$ and extreme points of the set of
distributions of \exch{} random infinite graphs
$H\in\cbloo$.
This correspondence is given by 
\begin{equation}
t(F,\gG)=\P(H\supset F)
\end{equation}
for every $F\in\BL$.
Furthermore, $H\restnn\pto \gG$ in $\cbq$ as $n_1,n_2\to\infty$.
\end{corollary}

\begin{remark}
  We have not checked whether $H\restnn\asto \gG$ in $\cbq$ as
  $n_1,n_2\to\infty$. This holds at least for a subsequence
  $(n_1(m),n_2(m))$ with both $n_1(m)$ and $n_2(m)$ non-decreasing
because then $\tinj(F,H\restnn)$ is a reverse martingale.
\end{remark}

\begin{theorem}
  \label{TE2B}
Let $H$ be an \exch{} random infinite bipartite graph.
Then the following are equivalent.
\begin{romenumerate}
  \item\label{te2Ba}
The distribution of $H$ is an extreme point in the set of \exch{}
distributions in $\cbloo$.
  \item\label{te2Bb'}
If $F_1$ and $F_2$ are two (finite) bipartite graphs with 
the vertex sets $V_j(F_1)$ and  $V_j(F_2)$ disjoint subsets of $\bbN$
for $j=1,2$, then
\begin{equation*}
  \P(H\supset F_1\cup F_2)=\P(H\supset F_1)\P(H\supset F_2).
\end{equation*}
\end{romenumerate}
\end{theorem}

The construction in \refS{Srep} takes the following form; note that
there is no need to assume symmetry of $W$.
For every $W\in\cw$, we define an infinite random bipartite graph
$\gwoooo\in\cbloo$ as follows: we
first choose two sequence $X_1,X_2,\dots$  and $Y_1,Y_2,\dots$ 
of \iid{} random variables
uniformly distributed on $\oi$, and then, given these sequences,
for each pair $(i,j)\in \bbN\times\bbN$ we draw an edge 
$ij$ with probability $W(X_i,Y_j)$, 
independently for
all pairs $(i,j)$.
Further, let $\gwnn$ be the restriction $\gwoooo\restnn$, which is
obtained by the same construction with finite sequences
$X_1,\dots,X_{n_1}$ and $Y_1,\dots,Y_{n_2}$.

It is evident that $\gwoooo$ is an \exch{} infinite random bipartite
graph. Furthermore, it satisfies \refT{TE2B}\ref{te2Bb'}.
\refT{TEB} and \refC{CEB} yield a corresponding element
$\uww\in\cboo\cong\cbxq$ such that $\gwnn\pto\uww$ as
$n_1,n_2\to\infty$ and, for every $F\in\BLkk$, 
\begin{equation}\label{tuww}
  t(F,\uww)
=
\int_{\oi^{k_1+k_2}} \prod_{ij\in\E(F)} W(x_i,y_j) 
  \dd x_1\dots \dd x_{k_1}  \dd y_1\dots \dd y_{k_2}.
\end{equation}

The result by Aldous \cite{Aldous} in the non-symmetric case
is that every \exch{} infinite random bipartite graph is obtained as a
mixture of such $\gwoooo$; in 
other words as $\gwoooo$ with a random $W$. 

By \refT{TEB} and \refC{CEB} above, this implies (and is implied by)
the fact that every element of $\cbq$ equals $\uww$ for some
(non-unique) $W\in\cw$;
the bipartite version of the characterization by \Lovaszetal{}
\cite{LSz}.

\section{Directed graphs}\label{Sdir}

A \emph{directed graph} $G$ consists of a vertex set
$V(G)$ and an edge set $E(G)\subseteq V(G)\times V(G)$; 
the edge indicators thus form an arbitrary zero--one matrix
\set{X_{ij}}, $i,j\in V(G)$. Note that we allow loops, corresponding
to the diagonal indicators $X_{ii}$.
The definitions and results above have analogues for directed
graphs too, with mainly notational differences. We sketch these in
this section, leaving the details to the reader.

Let $\DLn$ be the set of the $2^{n^2}$ labelled directed graphs
with vertex set $[n]$ and let $\DUn$ be the quotient
set $\DLn/\cong$ of unlabelled directed graphs with $n$
vertices; further, let
$\DL\=\bigcup_{n\ge1}\DLn$ and $\DU\=\bigcup_{n\ge1}\DUn$.

The definitions in \refS{Sdef} apply to directed graphs too, with at
most notational differences.
$G[k]$ and $G[k]'$ now are random directed graphs and $t(F,G)$,
$\tinj(F,G)$ and $\tind(F,G)$ are defined 
for (unlabelled) directed graphs $F$ and $G$ 
by \eqref{t}--\eqref{tind}.
We now define $\tau:\DU\to\oid$ by, \cf{} \eqref{tau}, 
\begin{equation}\label{taud}
  \tau(G)\=(t(F,G))_{F\in\cD}\in\oid. 
\end{equation}
We define $\cdx\=\tau(\cD)\subseteq\oid$ to be the image of $\cD$
under this mapping $\tau$, and let $\cdxq$ be the closure of $\cdx$ in
$\oid$; this is
a compact metric space.

Again, $\tau$ is not injective.
We let
$\cdd$ be the union of $\cD$ and some one-point set \set{*} and
consider the mapping $\tauu:\cD\to\oidd=\oid\times\oi$ defined by
\eqref{tau+} as before.
Then $\tauu$ is injective
and we can identify $\cD$ with its image $\tauu(\cD)\subseteq\oidd$
and define $\cdq\subseteq\oidd$ as its closure; this is a compact
metric space.
The functions $t(F,\cdot)$,
$\tinj(F,\cdot)$, $\tind(F,\cdot)$ and 
$v(\cdot)\qw$, for $F\in\cD$,
have unique continuous
extensions to $\cdq$.

We let $\cdoo\=\set{G\in\cdq:v(G)=\infty}$; 
this is the set of all limit objects of sequences
$(G_n)$ in $\cD$ with $v(G_n)\to\infty$.
By \eqref{a4xx}, $\tinj(F,G)=t(F,G)$ for every $G\in\cdoo$
and every $F\in\cD$. The projection $\pi:\cdq\to\cdxq$ restricts to a
homeomorphism $\cdoo\cong\cdxq$.

All results in Sections \ref{Sdef}--\ref{Sexch} are valid for directed
graphs too, with at most notational differences.

The main difference for the directed case concerns the representations
discussed in \refS{Srep}. Since two vertices may be connected by up to
two directed edges (in opposite directions), and the events that the
two possible edges occur typically are dependent, a single function
$W$ is no longer enough. Instead, we have a representation using 
several functions as follows.

Let $\WW$ be the set of quintuples
$\bW=(W_{00},W_{01},W_{10},W_{11},\wo)$ where  
$\wab:\oi^2\to\oi$ and $\wo:\oi\to\setoi$ are measurable functions 
such that $\sum_{\ga,\gb=0}^1\wab=1$ and $\wab(x,y)=W_{\gb\ga}(y,x)$
for $\ga,\gb\in\setoi$ and $x,y\in\oi$.
For $\bW\in\WW$, we define a random infinite directed graph
$\gbwoo$
by specifying its edge indicators $X_{ij}$
as follows: we
first choose a sequence $\xx_1,\xx_2,\dots$ of \iid{} random variables
uniformly distributed on $\oi$, and then, given this sequence,
let $X_{ii}=\wo(\xx_i)$ and 
for each pair $(i,j)$ with $i<j$ choose $X_{ij}$ and $X_{ji}$ at random
such that
\begin{equation}\label{dirw}
 \P(X_{ij}=\ga \text{ and } X_{ji}=\gb)=\wab(\xx_i,\xx_j),
\qquad \ga,\gb\in\setoi;
\end{equation}
this is done independently for
all pairs $(i,j)$ with $i<j$ (conditionally given \set{\xx_i}).
In other words, for every $i$ we draw a loop at $i$ if $\wo(\xx_i)=1$ 
and
for each pair $(i,j)$ with $i<j$ we 
draw edges $ij$ and $ji$ at random such that \eqref{dirw} holds.
Further, let $\gbwn$ be the restriction $\gbwoo\rest n$, which is
obtained by the same construction with a finite sequence $\xx_1,\dots,\xx_n$.

In particular, note that the loops appear independently, each with
probability $p=\P\bigpar{w(\xx_1)=1}$. We may specify the loops more clearly 
by the following alternative version of the construction.
Let $\csp\=\oi\times\setoi$ and
let $\WWW$ be the set of quadruples
$\bW=(W_{00},W_{01},W_{10},W_{11})$ where 
$\wab:\csp^2\to\oi$  are measurable functions 
such that $\sum_{\ga,\gb=0}^1\wab=1$ and $\wab(x,y)=W_{\gb\ga}(y,x)$
for $\ga,\gb\in\setoi$ and $x,y\in\csp$.
For every $\bW\in\WWW$ and $p\in\oi$, we define a
random infinite directed graph 
$\gbwpoo$
by specifying its edge indicators $X_{ij}$
as follows: We
first choose 
sequences $\xi_1,\xi_2,\dots$ and $\zeta_1,\zeta_2,\dots$ of random
variables, all independent, with $\xi_i\sim\uoi$ and
$\zeta_i\sim\Be(p)$, \ie, $\zeta_i\in\setoi$ with $\P(\zeta_i=1)=p$;
we let $Y_i\=(\xi_i,\zeta_i)\in\csp$.
Then, given these sequences,
let $X_{ii}=\zeta_i$ and 
for each pair $(i,j)$ with $i<j$ choose $X_{ij}$ and $X_{ji}$ at random
according to \eqref{dirw},
independently for
all pairs $(i,j)$ with $i<j$ (conditionally given \set{\xx_i}).
In other words, $\zeta_i$ is the indicator of a loop at $i$.
Further, let $\gbwpn$ be the restriction $\gbwpoo\rest n$, which is
obtained by the same construction with a finite sequence $\xx_1,\dots,\xx_n$.

It is obvious from the symmetry of the construction that the random
infinite directed graphs $\gbwoo$ and $\gbwpoo$ are \exch. Further,
using \refT{TE2}, their distributions are extreme points, so by
\refC{CE} they correspond to directed graph limits, \ie, elements of
$\cdoo$, which we denote by $\ubw$ and $\ubwp$, respectively;
\eqref{ce}
shows that if $F\in\cD_k$, then
\begin{align*}
  t(F,\ubw)=\P\bigpar{F\subseteq \gbwx k},
&&&
  t(F,\ubwp)=\P\bigpar{F\subseteq \gbwpx k}.
\end{align*}
By \refT{TE} and \refR{Ras}, $\gbwn\to\ubw$ and $\gbwpn\to\ubwp$ \as{}
as \ntoo. 

We can show a version of the representation
results in \refS{Srep} for directed graphs.

\begin{theorem}\label{TD}
An \exch{} random infinite directed graph is obtained as a mixture of
$\gbwoo$; in other words, as $\gbwoo$ with a random $\bW$. 
Alternatively, it is obtained as a mixture of 
$\gbwpoo$; in other words, as $\gbwpoo$ with a random $(\bW,p)$. 

Every directed graph limit, \ie, every element of $\cdoo$, is $\ubw$
for some $\bW\in\WW$, or equivalently $\ubwp$ for some $\bW\in\WWW$
and $p\in\oi$. 
\end{theorem}

\begin{proof}
For jointly \exch{} random arrays \set{X_{ij}} of zero--one variables, 
the Aldous--Hoover representation theorem takes the form 
\cite[Theorem 7.22]{Kallenberg:exch} 
  \begin{align*}
X_{ii}&=f_1(\xio,\xii), \\	
X_{ij}&=f_2(\xio,\xii,\xij,\xiij), \qquad i\neq j,	
  \end{align*}
where $f_1:\oi^2\to\setoi$ and $f_2:\oi^4\to\setoi$ are two measurable
functions, $\xi_{ji}=\xiij$, and 
$\xio$, $\xii$ ($1\le i$) and $\xiij$ ($1\le i<j$) are independent
random variables uniformly distributed on $\oi$ (as in the proof of
\refT{TU}). 
If further the distribution of the array \set{X_{ij}} is an extreme
point in the set of \exch{} distributions, then by \refT{TE2} and
\cite[Lemma 7.35]{Kallenberg:exch}, there exists such a representation
where $f_1$ and $f_2$ do not depend on $\xio$, so 
$X_{ii}=f_1(\xii)$ and $X_{ij}=f_2(\xii,\xij,\xiij)$, $i\neq j$.	
In this case, define $\wo=f_1$ and
\begin{equation*}
  \wab(x,y)\=\P\bigpar{f_2(x,y,\xi)=\ga\text{ and } f_2(y,x,\xi)=\gb},
\qquad \ga,\gb\in\setoi,
\end{equation*}
where $\xi\sim\uoi$. This defines a quintuple $\bW\in\WW$, such that
the edge indicators $X_{ij}$ of $\gbwoo$ have the desired distribution.

In general, the variable $\xio$ can be interpreted as making $\bW$
random.

To obtain the alternative representation, let
$\zeta_i\=w(\xi_i)=X_{ii}$ and $p\=\P(\zeta_i=1)$. There exists a
measure preserving map $\phi:(\csp,\mu_p)\to\oi$ such that
$\oi\times\set{j}$ is mapped onto $\set{x\in\oi:\wo(x)=j}$ for $j=0,1$
(\ie, $\wo\circ\phi(x,\zeta)=\zeta$), and we can use the quadruple 
$(\wab\circ\phi)_{\ga,\gb}$.

The representations for graph limits follow by \refC{CE} as discussed
above.
\end{proof}

\begin{example}
  A \emph{random tournament} $T_n$
is a random directed graph on $n$ vertices without loops where
each pair of vertices is connected by exctly one edge, with random
direction (with equal probabilities for the two directions, and
independent of all other edges). 
This equals $\gbwn$ or $\gbwpn$ with $W_{00}=W_{11}=0$,
$W_{01}=W_{10}=1/2$, and $w=0$ or $p=0$, and converges thus \as{} to
the limit $\gG_{\bW,0}$ for $\bW=(\wab)_{\ga,\gb}$.
\end{example}

Note that if \set{X_{ij}} are the edge indicators of an \exch{} random
infinite directed graph, then the loop indicators \set{X_{ii}} form a
binary \exch{} sequence, and the representation as $\gbwpoo$ in
\refT{TD} exhibits them as a mixture of \iid{} $\Be(p)$ variable,
which has brought us back to deFinetti's theorem \ref{deF}.

\newcommand\AAP{\emph{Adv. Appl. Probab.} }
\newcommand\JAP{\emph{J. Appl. Probab.} }
\newcommand\JAMS{\emph{J. \AMS} }
\newcommand\MAMS{\emph{Memoirs \AMS} }
\newcommand\PAMS{\emph{Proc. \AMS} }
\newcommand\TAMS{\emph{Trans. \AMS} }
\newcommand\AnnMS{\emph{Ann. Math. Statist.} }
\newcommand\AnnPr{\emph{Ann. Probab.} }
\newcommand\CPC{\emph{Combin. Probab. Comput.} }
\newcommand\JMAA{\emph{J. Math. Anal. Appl.} }
\newcommand\RSA{\emph{Random Struct. Alg.} }
\newcommand\ZW{\emph{Z. Wahrsch. Verw. Gebiete} }
\newcommand\DMTCS{\jour{Discr. Math. Theor. Comput. Sci.} }

\newcommand\AMS{Amer. Math. Soc.}
\newcommand\Springer{Springer-Verlag}
\newcommand\Wiley{Wiley}

\newcommand\vol{\textbf}
\newcommand\jour{\emph}
\newcommand\book{\emph}
\newcommand\inbook{In \emph}
\def\no#1#2,{\unskip#2, no. #1,} 
\newcommand\toappear{\unskip, to appear}

\newcommand\webcite[1]{%
   \penalty1\texttt{\def~{{\tiny$\sim$}}#1}\hfill\hfill}
\newcommand\webcitesvante{\webcite{http://www.math.uu.se/~svante/papers/}}
\newcommand\arxiv[1]{\webcite{http://arxiv.org/#1}}

\def\nobibitem#1\par{}

\end{document}